\documentclass[journal]{IEEEtran}
\usepackage{cite}

\ifCLASSINFOpdf
\else
\fi
\usepackage{amsmath,graphicx}
\usepackage{hyperref}
\usepackage{graphicx}
\usepackage{float}
\usepackage{amsfonts}
\usepackage{amsmath,amssymb}
\usepackage{epstopdf}
\usepackage[ruled]{algorithm2e}
\usepackage{setspace}
\usepackage{color}
\usepackage{wrapfig}
\usepackage{marvosym}
\usepackage{cleveref}
\usepackage{ragged2e}
\usepackage{cite}
\usepackage[scaled]{helvet} 
\usepackage{fourier}
\usepackage{pdflscape}
\usepackage[english]{babel}
\usepackage{multirow}
\usepackage[usenames,dvipsnames,svgnames,table]{xcolor}
\usepackage{nomencl}
\usepackage{enumitem}

\usepackage{algorithmic}

\makenomenclature

\newcommand{\bi}{\boldsymbol{i}}

\newcommand{\fd}{\mathfrak{d}}
\newcommand{\fD}{\mathfrak{D}}

\newtheorem{thm}{Theorem}[section]

\begin{document}
%
\title{A Two-Stage Decomposition Approach for AC Optimal Power Flow}
%
%
%

\author{Shenyinying~Tu, ~Andreas W\"achter, ~Ermin~Wei 
\thanks{Manuscript submitted December 18, 2019. This work was supported by Leslie and Mac McQuown.   The second author was partially supported by National Science Foundation grant DMS-1522747 and by Los Alamos National Laboratory, its Center for Nonlinear Studies, and its Ulam Scholarship program. Los Alamos National Laboratory is operated by Triad National Security, LLC, for the National Nuclear Security Administration of U.S. Department of Energy (Contract No. 89233218CNA000001) under the auspices of the NNSA of the U.S. DOE at LANL under Contract No. DE- AC52- 06NA25396..}
\thanks{The authors are with Northwestern University, Evanston, IL 60201 USA (email: shenyinyingtu2021@u.northwestern.edu; waechter@iems.northwestern.edu; ermin.wei@northwestern.edu)}}

\IEEEpubid{}
\maketitle
\begin{abstract}
The alternating current optimal power flow (AC-OPF) problem is critical to power system operations and planning, but it is generally hard to solve due to its nonconvex and large-scale nature. This paper proposes a scalable decomposition approach in which the power network is decomposed into a master network and a number of subnetworks, where each network has its own AC-OPF subproblem. This formulates a two-stage optimization problem and requires only a small amount of communication between the master and subnetworks. The key contribution is a smoothing technique that renders the response of a subnetwork differentiable with respect to the input from the master problem, utilizing properties of the barrier problem formulation that naturally arises when subproblems are solved by a primal-dual interior-point algorithm.  Consequently, existing efficient nonlinear programming solvers can be used for both the master problem and the subproblems. The advantage of this framework is that speedup can be obtained by processing the subnetworks in parallel, and it has convergence guarantees under reasonable assumptions. The formulation is readily extended to instances with stochastic subnetwork loads. Numerical results show favorable performance and illustrate the scalability of the algorithm which is able to solve instances with more than 11 million buses. 
\end{abstract}


%
\IEEEpeerreviewmaketitle

\section{Introduction}
%
%
%
%
\IEEEPARstart{T}{he} study of optimal power flow (OPF) \cite{huneault1991survey,frank2012optimal,capitanescu2011state,molzahn2019survey} is fundamental in power systems, because it is an essential building block to investigate questions in 
 operation and planning, such as unit commitment \cite{sen1998optimal,padhy2004unit}, stability and reliability assessment \cite{canizares2002voltage,sauer1998power}, etc. It seeks to optimize some cost function, such as generation cost or transmission loss, while also satisfying the physical constraints of the power network. 

One of the difficulties of solving alternating current optimal power flow (AC-OPF) problems arises from the non-convex constraint set. The direct current optimal power flow (DC-OPF) model has been widely used in practice by linearizing the AC-OPF model \cite{alsac1990further, coffrin2014linear}, which assumes the optimal solution is close to the normal operating point. However, the penetration of renewable energy introduces high fluctuations in the energy generation, which can deviate significantly from the regular operation. Hence the DC-OPF formulation with renewable generation may lead to infeasible or sub-optimal solutions. 
On the other hand, the growing network size makes AC-OPF problems computationally expensive to solve. To efficiently implement OPF in large-scale systems, it is beneficial to decompose the overall problem into smaller pieces, each of which can be solved independently. 

One application of this decomposition approach is the co-optimization of a transmission network and distribution networks. Currently, the transmission system operator (TSO) and the distribution system operator (DSO) administer their networks independently with little coordination \cite{li2016coordinated}. In a traditional distribution system, power flow is unidirectional, and distribution systems are typically modeled as a load bus in the transmission system. Similarly, a DSO operates the distribution system by simplifying the transmission system as a voltage source.
However, with high penetration of distributed energy resources (DERs), it is also reasonable to consider each distribution system as an active power plant \cite{jia2015hierarchical}.  This creates the need to operate the transmission system and the distribution systems in a more coordinated way. The current interactions between TSO and DSO are reviewed in \cite{zegers2014tso}. The authors concluded that bidirectional communication between TSO and DSO is needed. In this case, the distribution system can no longer be decoupled from the transmission network, resulting in a challenging large-scale AC-OPF problem.

Various decomposition approaches for power networks have been proposed. By applying a decomposition technique, a large network is divided into smaller subnetworks that can be solved efficiently and in parallel \cite{molzahn2017survey}. Decentralized OPF approaches are proposed in \cite{kim1997coarse,baldick1999fast, kim2000comparison, nogales2003decomposition, hug2009decentralized}, where augmented Lagrangian methods are employed, including the auxiliary problem principle, the alternating direction multiplier method (ADMM) and the predictor-corrector proximal multiplier method. Each region solves its own OPF subproblem independently and communicates with its neighbors defined by respective partitions. However, augmented Lagrangian algorithms may fail to converge due to the lack of convexity  \cite{erseghe2014distributed}.

Primal-dual interior point methods (PDIPM) for AC-OPF problems are studied in \cite{jabr2002primal,capitanescu2007interior}. However, these methods require the computation of a Newton step from a large linear system, which might be prohibitively expensive in large-scale OPF problems. One recent work \cite{minot2016parallel} proposes a parallel PDIPM for decomposed power networks based on matrix splitting \cite{varga2009matrix}, in which an outer PDIPM loop and an inner matrix-splitting loop are involved. Though the paper shows that the problem can converge with a few PDIPM iterations, it still takes hundreds of matrix-splitting iterations within each PDIPM iteration.

More recent works on distributed algorithms \cite{dall2013distributed,lam2012distributed,kraning2013dynamic} focus on relaxed models of OPF problems, where the nonconvex constraint set is relaxed by a convex outer approximation. Although these algorithms provide improvements in the computation time, feasibility is not always guaranteed.

Most of the existing studies on decomposition algorithms consider multi-area transmission systems \cite{molzahn2017survey}, and most work on the optimization of distribution systems does not consider the interaction with the transmission system \cite{low2013convex}.  The studies that examine decomposition algorithms for TSO-DSO interactions include a master-slave iterative algorithm \cite{sun2014master,li2016coordinated2}, which decompose the optimality conditions of the entire problem into those for the transmission and distribution networks and propose a heuristic to alternate between solving the two subproblems. Furthermore, hierarchical coordination is studied in \cite{yuan2017hierarchical}, which applies Benders decomposition.
Because the AC power flow equations in the distribution systems are convexified using second-order cone constraints, the solutions obtained are often not feasible for the original AC problem. Finally, a Lagrangian relaxation method is extended to TSO-DSO cooperation problems in \cite{mohammadi2018diagonal}, where the Lagrangian term is linearized to enable parallel computations for both transmission and distribution networks, but it requires an extra outer loop for the Lagrangian term to converge. Among those, only \cite{yuan2017hierarchical} offers theoretical convergence guarantees.  The largest network considered by any of these publications corresponds to an undecomposed system with 1,088 buses \cite{li2016coordinated2}, which is several orders of magnitudes smaller than the 11 million bus system solved here.

In this paper, we propose a novel way of decomposing power networks and solve it using a new two-stage optimization technique. The approach is inspired by \cite{demiguel2006local,demiguel2008decomposition}, which considers only local convergence properties. Here, we extend this idea to a practical algorithm with global convergence guarantees and assess its performance in a power flow setting. In particular, the network is partitioned into a master network and a set of subnetworks, each having its own AC-OPF problem. At each iteration of the master problem, subproblems are solved in parallel using an PDIPM, and the master network then makes decisions based on the output of subnetworks. Communication between the master network and subnetworks is only required for a small amount of variables. 

In general, the optimal response of a subproblem is non-differentiable with respect to the communication variables.  \textit{Our key contribution is a smoothing technique that results in a differentiable optimal value function of a subproblem as a function of the input from master problem. This permits the use of efficient general-purpose gradient based nonlinear programming (NLP) solvers for the solution of the master problem.  This is made possible by utilizing properties of the barrier problem formulation that naturally arises when the subnetworks are solved with an PDIPM.}

We apply this decomposition approach to AC-OPF problems with TSO-DSO interactions, by considering the transmission network as the master problem and distribution networks as subnetworks. We follow the current assumption in literature \cite{sun2014master,li2016coordinated2,yuan2017hierarchical,mohammadi2018diagonal} where TSO and DSO are cooperative and collectively optimizing the overall system performance. The distribution networks are not necessarily tree structured networks. Our framework is applicable to the case with meshed distribution networks.

Moreover, because of the uncertainties of DERs power production in distribution systems \cite{su2012energy}, the data in subproblems is often stochastic. 
To account for this, the proposed framework can naturally approximate the expected cost by making many copies of a subnetwork, each one with different realizations of the random data. All of the sampled distribution network problems can be solved in parallel, which makes the resulting very large-scale formulation possible to solve.

This paper is organized as follows: Section \ref{sec:ACOPF} reviews the branch flow model formulation of AC-OPF problems. Section \ref{sec:bilevel} proposes the decomposition approach of the AC-OPF problem and reformulates the problem as a two-stage optimization problem. Implementation details and numerical results are given in Section \ref{sec:numerical}, and Section \ref{sec:conclusion} concludes the paper.

\section{Branch flow model} \label{sec:ACOPF}
Consider a directed power network $G := (N,E)$, where $N$ denotes the set of buses, and $E$ denotes the set of branches. For each branch $(i,j)$ in $E$, let $y_{ij}$ denote its admittance, $z_{ij}:= 1/y_{ij}$ be the corresponding impedance. Let $I_{ij}$ be the complex current, and $S_{ij}: = P_{ij} + \bi Q_{ij}$ be the complex power from bus $i$ to bus $j$. For each bus $j$ in $N$, let $V_j$ be the complex voltage, and $s_j$ be the power injection. Then the branch flow model is defined by \cite{low2013convex}:
\begin{subequations}\label{eq:original}
\begin{align}
& I_{ij} = y_{ij}(V_i - V_j), ~~\forall i\to j \in E \label{eq:1_current} \\
& S_{ij} = V_i I_{ij}^H, ~~\forall i\to j \in E \label{eq:1_powerFlow} \\
& s_j = \sum_{k:j\to k} S_{jk} - \sum_{i:i\to j} (S_{ij} - z_{ij}|I_{ij}|^2), ~~\forall j\in N \label{eq:1_powerInjection} \\
& \underline{v}_j \leq |V_j|^2\leq \overline{v}_j, ~~\forall j\in N \label{eq:1_voltage_limit} \\
& \underline{s}_j \leq s_j \leq \overline{s}_j, ~~\forall j\in N \label{eq:1_power_limit}
\end{align}
\end{subequations}
where the complex conjugate of $I_{ij}$ is denoted by $I^H_{ij}$. 
At node $j$, $\underline{v}_j$, $\overline{v}_j$ give the lower and upper bounds on voltage magnitude, and $\underline{s}_j$, $\overline{s}_j$ give the lower and upper bounds on power injection. Typical costs of the OPF problem include the line loss along the branches and/or power generation cost, which are usually quadratic functions of the variables. Let $C(V, S, I)$ be the cost function of the problem. Then we can formulate the AC-OPF problem as 
\begin{equation} \label{eq:meshACOPF}
\begin{split}
\min_{V, S, I, s}  & \quad C(V, S, I, s) \\
s.t. & \quad \eqref{eq:original}.
\end{split}
\end{equation}

The rest of the paper will focus on solving this problem. Note that our decomposition approach is not restricted to the branch flow model, and it can be applied to other power flow model representations.

\section{Two-stage Optimization Framework} \label{sec:bilevel}

\subsection{Decomposition Scheme} \label{sec:decomposition}
\begin{figure}[!t]
  \centering
  \includegraphics[width=0.34\textwidth]{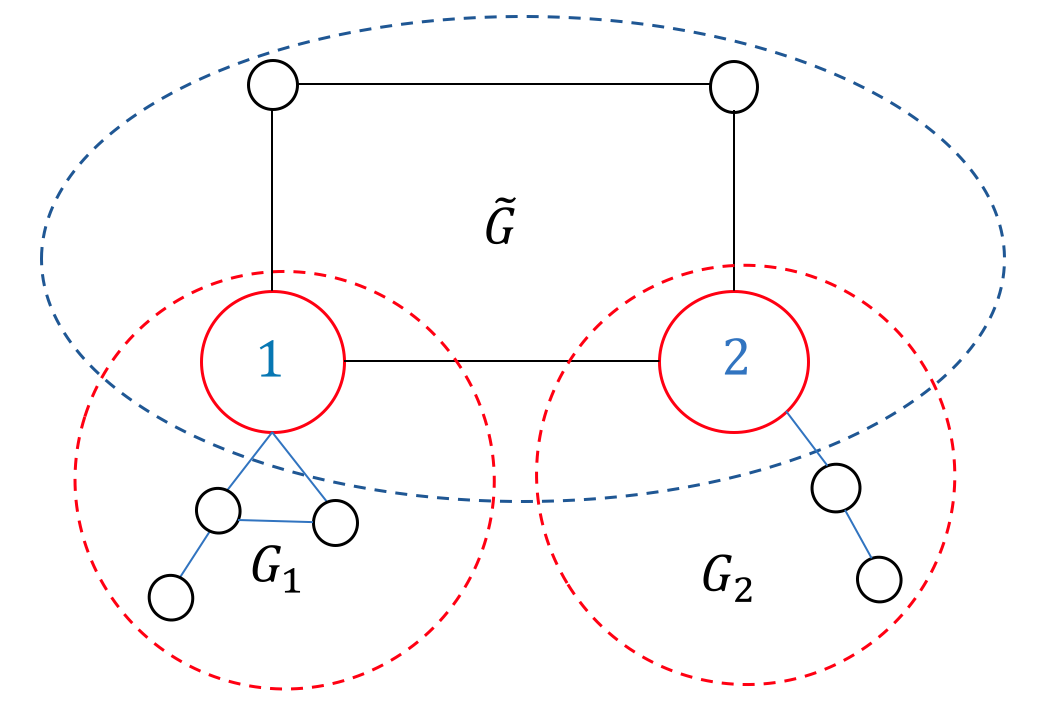}
  \caption{Original network}
  \label{fig:undecomposed}
\end{figure}

\begin{figure}[!t]
  \centering
  \includegraphics[width=0.33\textwidth]{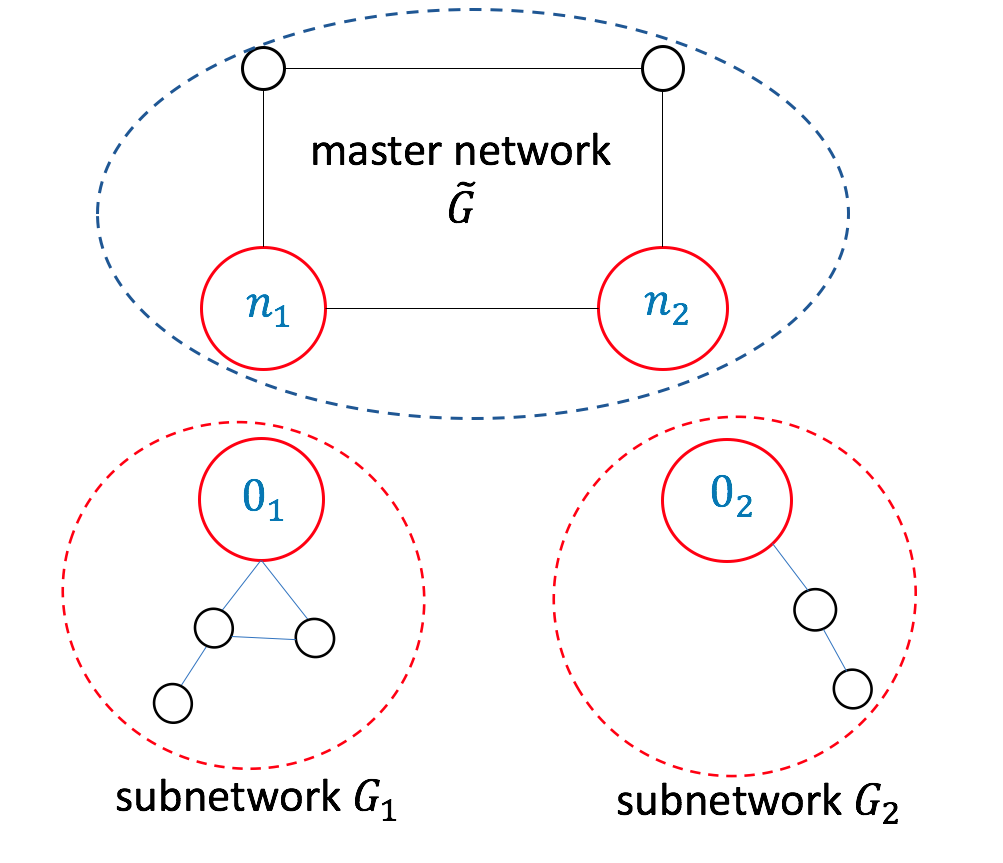}
  \caption{Decomposed network}
  \label{fig:decomposed}
\end{figure}

This paper proposes a decomposition approach of solving the AC-OPF problem \eqref{eq:meshACOPF}. In our decomposition scheme, a power network is decomposed into a master network and a number of independent subnetworks $\{\fd: \fd \in \fD\}$, as sketched in Fig~\ref{fig:undecomposed}. 
 Let $\tilde{G} := (\tilde{N}, \tilde{E})$ be the master network, and let $G_{\fd}: = (N_\fd, E_\fd)$ represent the subnetwork, for $\fd \in \fD$. Motivated by the topology of transmission and distribution networks, we suppose that the network can be decomposed such that each subnetwork $G_\fd$ overlaps with the master network by exactly one bus $\fd$, denoted by $0_\fd \in N_\fd$ or $n_\fd \in \tilde{N}$ as shown in Fig~\ref{fig:decomposed}. To avoid confusion in notation, we mark all master problem variables with a tilde (`` $\tilde{}$ '').
In principle, the proposed approach can be applied when more than one bus is shared, see Section~\ref{sec:penalty}.

After decomposition, for each connecting bus, the problem has decision variables for both the master network $(\tilde{V}_{n_\fd}, \tilde{s}_{n_\fd})$ and the corresponding subnetwork $(V_{0_\fd}, s_{0_\fd})$. For the master network, let $\tilde{V}_{n_\fd}$ denote the complex voltage at $n_\fd$, and $\tilde{s}_{n_\fd}$ be the net power flow from the master network into bus $n_\fd$:
\begin{equation}
\tilde{s}_{n_\fd} = \sum_{n_\fd\to k \in \tilde{E}} \tilde{S}_{n_\fd k} - \Big(\sum_{i\to n_\fd \in \tilde{E}}(\tilde{S}_{i n_\fd}-\tilde{z}_{i n_\fd}|\tilde{I}_{i n_\fd}|^2)\Big).
\end{equation}
\indent For subnetwork $\fd$, let $V_{0_\fd}$ denote the complex voltage at $0_\fd$, and $s_{0_\fd}$ be the net power flow from subnetwork $\fd$ into bus $0_\fd$:
\begin{equation}
s_{0_\fd} = \sum_{0_\fd\to k \in E_\fd} S_{0_\fd k} - \Big(\sum_{i\to 0_\fd \in E_\fd}(S_{i 0_\fd}-z_{i 0_\fd}|I_{i 0_\fd}|^2)\Big).
\end{equation}
\indent To be consistent with the undecomposed formulation \eqref{eq:original}, we impose boundary conditions that couple the decision variables at the connecting bus $n_\fd$:
\begin{subequations}
\begin{align}
& \tilde{V}_{n_\fd}  = V_{0_\fd}, \\
& \tilde{s}_{n_\fd} + s_{0_\fd} = 0.
\end{align}
\end{subequations}
\indent Therefore, the entire AC-OPF problem can be decomposed into a master level problem and a set of subproblems:
\begin{equation*} 
\begin{split}
\min &\quad C(\tilde{V},\tilde{S},\tilde{I},\tilde{s}) + \sum_{\fd \in \fD} C^*_\fd(\tilde{V}_{n_\fd}, \tilde{s}_{n_\fd}) \\
s.t. & \quad \tilde{I}_{ij} = \tilde{y}_{ij}(\tilde{V}_i - \tilde{V}_j), ~~\forall i\to j \in \tilde{E} \\
&  \quad \tilde{S}_{ij} = \tilde{V}_i \tilde{I}_{ij}^H, ~~\forall i\to j \in \tilde{E} \\
&  \quad \tilde{s}_j = \sum_{j\to k} \tilde{S}_{jk} - \sum_{i\to j} (\tilde{S}_{ij} - \tilde{z}_{ij}|\tilde{I}_{ij}|^2), ~~\forall j\in \tilde{N} \setminus \{n_\fd \}_{\fd \in \fD}  \\
& \quad \tilde{s}_{n_\fd} = \sum_{n_\fd\to k \in \tilde{E}} \tilde{S}_{n_\fd k} - \Big(\sum_{i\to n_\fd \in \tilde{E}}(\tilde{S}_{i n_\fd}-\tilde{z}_{i n_\fd}|\tilde{I}_{i n_\fd}|^2)\Big)\\
&  \quad \underline{\tilde{v}}_j \leq |\tilde{V}_j|^2\leq \overline{\tilde{v}}_j, ~~\forall j\in \tilde{N} \\
&  \quad \underline{\tilde{s}}_j \leq \tilde{s}_j \leq \overline{\tilde{s}}_j, ~~\forall j\in \tilde{N} \setminus \{n_\fd \}_{\fd \in \fD} \\
\end{split}
\end{equation*}
where $C^*_\fd(\tilde{V}_{n_\fd}, \tilde{s}_{n_\fd})$ is the optimal objective value of the following subproblem $\fd$ given fixed values for $\tilde{V}_{n_\fd}$, $\tilde{s}_{n_\fd}$:
\begin{subequations} \label{eq:subproblem_original}
\begin{align}
\min & \quad C_\fd(V,S,I,s) \\
s.t & \quad V_{0_\fd} = \tilde{V}_{n_\fd}, \label{eq:constraint_v_equal} \\
& \quad s_{0_\fd} = -\tilde{s}_{n_\fd}, \label{eq:constraint_s_equal} \\ 
& \quad I_{ij} = y_{ij}(V_i - V_j), ~~\forall i\to j \in E_\fd \label{eq:constraint_IYV} \\
& \quad S_{ij} = V_i I_{ij}^H, ~~\forall i\to j \in E_\fd \label{eq:constraint_SVI} \\
& \quad s_j = \sum_{j\to k} S_{jk} - \sum_{i\to j} (S_{ij} - z_{ij}|I_{ij}|^2), ~~\forall j\in N_\fd \setminus 0_\fd \label{eq:constraint_sj} \\
& \quad s_{0_\fd} = \sum_{0_\fd\to k \in E_\fd} S_{0_\fd k} - \sum_{i\to 0_\fd \in E_\fd}(S_{i 0_\fd}-z_{i 0_\fd}|I_{i 0_\fd}|^2) \label{eq:constraint_s0} \\
& \quad \underline{v}_j \leq |V_j|^2\leq \overline{v}_j, ~~\forall j\in N_\fd \label{eq:constraint_v_bound} \\
& \quad \underline{s}_j \leq s_j \leq \overline{s}_j, ~~\forall j\in N_\fd \setminus 0_\fd \label{eq:constraint_s_bound}
\end{align}
\end{subequations}
\indent For simplification, let $x: = (\tilde{V}_i, ~\tilde{s}_i: i \in \tilde{N}; ~ \tilde{S}_{ij},~ \tilde{I}_{ij}: (i,j) \in \tilde{E})$ consists of all the variables in the master network. Let $x_\fd: = (\tilde{V}_{n_\fd}, \tilde{s}_{n_\fd})$ be the subset of master problem variables $x$ that couple subnetwork $\fd$ with the master network, and $y_\fd: = (V_i, ~s_i: i \in N_\fd; ~ S_{ij},~I_{ij}: (i,j) \in E_\fd)$ be the local variables in subnetwork $\fd$. Then the above AC-OPF problem \eqref{eq:upper} can be expressed compactly as a two-stage nonlinear programming problem:
\begin{subequations} \label{eq:upper}
\begin{align}
\min_{x} &\quad  C(x) + \sum_{\fd \in \fD} C^*_\fd(x_\fd) \label{eq:upper_obj}\\
s.t. & \quad g(x) = 0, \\
& \quad h(x) \leq 0,
\end{align}
\end{subequations}
where
\begin{subequations} \label{eq:subproblem}
\begin{align}
C^*_\fd(x_\fd) = \min_{y_\fd} & \quad  C_\fd(y_\fd;x_\fd) \\
s.t. & \quad g_\fd(y_\fd; x_\fd) = 0,\label{eq:subproblem_g} \\
& \quad  h_\fd(y_\fd;x_\fd) \leq 0. \label{eq:subproblem_const}
\end{align}
\end{subequations}
Here, the constraint functions $g$, $h$, $g_\fd$ and $h_\fd$ are smooth.

Even though, in the above formulation, only voltage magnitudes and power injections of the connecting bus are passed to a subnetwork, choices of communication variables can be made differently, such as voltage phase angle and current injection.
In principle, the proposed decomposition algorithm is applicable to any two-stage problem \eqref{eq:upper}--\eqref{eq:subproblem} with continuous variables in both stages.  This includes  other variations of the AC-OPF problem, such as subnetworks with stochastic wind turbines \cite{ahmadi2013multi} and three phase representation  \cite{molzahn2017survey}.  It also captures the decomposition of networks with more than one coupling bus, see Section~\ref{sec:penalty}.
These are settings that can be explored in future work.

%

\begin{figure}[!t]
  \centering
  \includegraphics[width=0.4\textwidth]{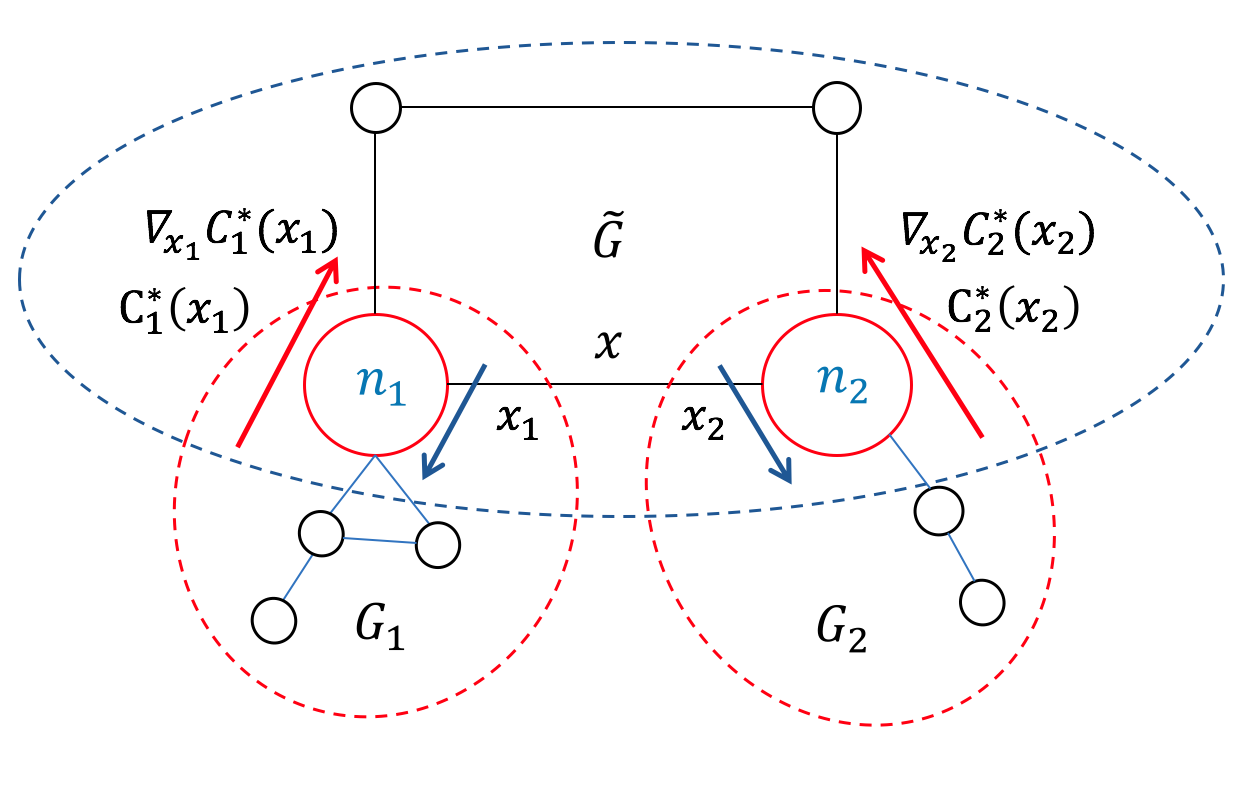}
  \caption{Network decomposition}
  \label{fig:bilevel}
\end{figure}

In the proposed method, the master problem \eqref{eq:upper} is optimized with a gradient-based nonlinear programming (NLP) solver. Whenever the NLP solver requires the value or derivatives of the objective function \eqref{eq:upper_obj}, the quantities $\{x_\fd\}_{\fd \in \fD}$ corresponding to the current iterate are passed to the subproblems \eqref{eq:subproblem}. Their optimal solutions are computed, and the optimal objective values $C^*_\fd(x_\fd)$, together with the first and second order derivatives $\frac{\partial C^*_\fd}{\partial x_\fd}$ and $\frac{\partial^2 C^*_\fd}{\partial x_\fd^2}$, are passed back to the NLP solver, which continues to solve \eqref{eq:upper} until a minimizer is found. The information exchange is illustrated in Fig \ref{fig:bilevel}. Since only few variables are communicated within the network, this approach can be efficiently implemented in a distributed setting. Note that this algorithm is ill-defined if any of the subproblems becomes infeasible given $x_\fd$ at one of the master problem iterates. For now, we assume subproblems are always feasible for any $x_\fd$. We will discuss in Section \ref{sec:penalty} how this assumption can be lifted.
\vspace{-15pt}
\subsection{Subproblems} \label{sec:subproblem}
\subsubsection{Smoothing of Subproblems} \label{sec:subproblem_smooth}
The description above ignores the crucial fact that, in general, the optimal value functions $C_\fd^*$ is not differentiable at some values of $x_\fd$.  More specifically, whenever the set of inequality constraints \eqref{eq:subproblem_const} that are tight at the optimal solution changes with $x_\fd$, the function $C_\fd^*$ typically is non-differentiable and experiences abrupt changes in first derivatives around those positions. This may result in convergence failures of the master problem NLP solver.

As a remedy, we replace \eqref{eq:subproblem} by its barrier problem formulation: 

\begin{subequations} \label{eq:relaxed_subproblems}
\begin{align}
C^*_\fd(x_\fd, \tilde{\mu}) = \min_{y^\fd} \quad & C_\fd(y_\fd; x_\fd) - \tilde{\mu}\sum_i \ln (s_i) \label{eq:relaxed_subproblems_obj} \\
s.t. \quad & h_\fd(y_\fd; x_\fd) +s = 0. \label{eq:ineq_constraint}
\end{align}
\end{subequations}
Here, the inequality constraint \eqref{eq:subproblem_const} has been replaced by an equality constraint that introduces slack variables $s$. (To simplify notation, we drop the  equality constraints \eqref{eq:subproblem_g} for the remainder of this section).
The objective function \eqref{eq:relaxed_subproblems_obj} now includes a logarithmic barrier term with weight $\tilde\mu>0$ which keeps the slack variables strictly positive.  
Interior point methods are based on this formulation and obtain an optimal solution of the original problem \eqref{eq:subproblem} by solving a sequence of barrier problems \eqref{eq:relaxed_subproblems} in which the barrier parameter $\tilde \mu$ converges to zero, thus recovering the solution to problem \eqref{eq:subproblem}, see Section~\ref{eq:2stagealg}.

The crucial property in our context is that \eqref{eq:relaxed_subproblems} no longer contains inequality constraints and that the optimal value function $C^*_\fd(x_\fd, \tilde{\mu})$ depends smoothly on $x_\fd$ (see Section~\ref{sec:2stagecom}).  This allows us a to define a smooth modified master problem,
\begin{equation} \label{eq:relaxed_master}
  \begin{split}
    \min_{x} & \quad C(x) + \sum_{\fd \in \fD} C^*_\fd(x_\fd, \tilde{\mu}) \\
    s.t. & \quad g(x) = 0 \\
    & \quad h(x) \leq 0,\\
  \end{split}
\end{equation}
which can be optimized with standard NLP solvers.  Section \ref{sec:2stagecom} describes how the derivatives of $C^*_\fd(x_\fd, \tilde{\mu})$ can be computed.

\subsubsection{Primal-Dual Interior Point Method}\label{sec:PDIP}
Now that the subproblem has the form of a barrier problem \eqref{eq:relaxed_subproblems}, it is natural to apply a PDIPM to solve it.  These methods work with the primal-dual optimality conditions for \eqref{eq:relaxed_subproblems},
\begin{equation}\label{eq:KKT}
 F(y_\fd, s_\fd,\lambda_\fd;x_\fd,\mu) = 0,
\end{equation}
where
$$F(y_\fd,s_\fd,\lambda_\fd;x_\fd, \mu) = 
\begin{pmatrix}
  \nabla_{y_\fd}  C_\fd(y_\fd;x_\fd)- \nabla_{y_\fd} h_{\fd}(y_\fd;x_\fd)^T \lambda \\
  h_\fd(y_\fd;x_\fd) + s_\fd \\
  s_\fd \circ \lambda_\fd - \mu \mathbf{1}
\end{pmatrix}.$$
Here, $\lambda$ is the dual variable corresponding to the constraint \eqref{eq:ineq_constraint}, $s \circ \lambda$ is the component-wise product of two vectors $s$ and $\lambda$, and the notation $\mathbf{1}$ defines a column vector with entries 1.
We can use any general purpose PDIPM solver, which seeks to find a root of $F(\,\cdot\,;x_\fd, \mu)$ for a fixed value of $\mu$ by applying Newton's method to \eqref{eq:KKT} and decreases $\mu$ to zero in order to converge to a local optimum (or at least a stationary point) of the original problem \eqref{eq:subproblem}. 

\begin{algorithm}
  \caption{Generic PDIPM Framework}\label{Alg:inner_alg}
  \textbf{Input:}  $x_\fd$, initial iterate $(y_{\fd}^{(0)}, s_\fd^{(0)},\lambda_\fd^{(0)})$, initial barrier parameter $\mu^{(0)}$.
  \begin{algorithmic}[1]
    \STATE Set $k\gets 0$.
    \REPEAT
      \WHILE{$\|F(y_\fd^{(k)}, s_\fd^{(k)},\lambda_\fd^{(k)};x_\fd,\mu^{(k)})\|>\epsilon_F^{(k)}$}
      \STATE Compute Newton step $(\Delta y_{\fd}^{(k)}, \Delta s_\fd^{(k)},\Delta \lambda_\fd^{(k)})$.\label{s:Newton}
      \STATE Perform line search to compute step size $\alpha^{(k)}\in(0,1]$.
      \STATE Update iterate
      $
      (y_{\fd}^{(k+1)}, s_\fd^{(k+1)},\lambda_\fd^{(k+1)}) = (y_{\fd}^{(k)}, s^{(k)}_\fd,\lambda^{(k)}_\fd) + \alpha^{(k)}(\Delta y_{\fd}^{(k)}, \Delta s^{(k)}_\fd,\Delta \lambda_\fd^{(k)}).
      $
      \STATE Increase iteration counter $k\gets k+1$.
      \ENDWHILE
      \STATE Decrease $\mu^{(k)}$. \label{s:muupdate} 
    \UNTIL{$\mu^{(k)}\leq\epsilon_{\mu}$.}
    \RETURN{$(y_{\fd}^{(k)}, s_\fd^{(k)},\lambda_\fd^{(k)})$.}
  \end{algorithmic}
\end{algorithm}

Algorithm \ref{Alg:inner_alg} shows the steps of a generic line search PDIPM.  Practical methods are more involved \cite{wachter2006implementation,waltz2006interior}, but this description highlights the features that are relevant in our context.  In the while loop, $\epsilon_F^{(k)}$ is the tolerance to which the barrier problem for the current value of $\mu^{(k)}$ is solved, and $\epsilon_{\mu}$ is the overall convergence tolerance. In a regular setting, $\epsilon_{\mu}$ is set to a tight tolerance $\epsilon$ (e.g., $10^{-8}$) and $\epsilon_F^{(k)}=0.1\mu^{(k)}$ \cite{wachter2006implementation}.

Step~\ref{s:Newton} requires the computation of the Newton step, which is computed as the solution of the linear system
\begin{equation} \label{eq:Newton}
J(y_{\fd}^{(k)}, s_\fd^{(k)}, \lambda_\fd^{(k)}; x_\fd, \mu^{(k)})
\begin{bmatrix}
\Delta y_\fd^{(k)} \\ \Delta s_\fd^{(k)} \\ \Delta \lambda_\fd^{(k)}
\end{bmatrix} = - F(y_\fd^{(k)}, s_\fd^{(k)},\lambda_\fd^{(k)};x_\fd,\mu^{(k)}) 
\end{equation}
where
\begin{equation} \label{eq:Jacobian}
J = 
\setlength\arraycolsep{2pt}
\begin{bmatrix} 
\nabla_{y_\fd}^2 C_\fd - \sum_i \nabla_{y_\fd}^2 h_{i_\fd} \lambda_i & 0 & -\nabla_{y_\fd} h^T \\
\nabla_{y_\fd} h             & I & 0                \\
  0       & \text{diag} (\lambda_\fd) & \text{diag}(s_\fd)
\end{bmatrix} 
\end{equation}
is the Jacobian of $F$. Here we drop the function arguments for brevity.

As we discussed in Section \ref{sec:subproblem_smooth}, the smoothed master problem \eqref{eq:relaxed_master} requires the solution of the barrier problem \eqref{eq:relaxed_subproblems} for a fixed given value $\tilde\mu>0$ of the barrier parameter.  This can be computed with a simple modification of Algorithm~\ref{Alg:inner_alg}.  Instead of decreasing the barrier parameter $\mu^{(k)}$ all the way to zero, the algorithm eventually fixes it to $\tilde\mu$ in Step~\ref{s:muupdate}, and from then on sets $\epsilon_F^{(k)}$ to the overall tight convergence tolerance $\epsilon$.  The solution returned is then an optimal primal-dual solution for \eqref{eq:relaxed_subproblems}.

We emphasize that this approach requires only a small modification of the termination criteria of a PDIPM. 
This has the significant practical advantage that powerful and efficient implementations of PDIPM such as Ipopt \cite{wachter2006implementation} and KNITRO \cite{waltz2006interior} could easily be adapted and utilized for the solution of the AC-OPF problem \eqref{eq:relaxed_subproblems}.  The next section shows that also the computation of derivatives of $C^*_\fd(x_\fd, \tilde{\mu})$ can exploit existing features of a PDIPM implementation.

\subsubsection{Derivative Computations}\label{sec:2stagecom}

Recall that $J$ is the Jacobian \eqref{eq:Jacobian} of primal-dual optimality conditions of problem \eqref{eq:relaxed_subproblems} for any given $x_\fd$ and $\tilde{\mu}$. We assume that some standard second-order optimality conditions, which typically hold for non-degenerate AC-OPF problems {\cite{hauswirth2018generic}}, are satisfied for \eqref{eq:relaxed_subproblems} so that $J$ is nonsingular. By the implicit function theorem, there exists a set of unique differentiable functions $(y_\fd^*(x_\fd), s_\fd^*(x_\fd), \lambda_\fd^*(x_\fd))$ in the neighborhood of $x_\fd$, where $y_\fd^*(x_\fd), s_\fd^*(x_\fd), \lambda_\fd^*(x_\fd)$ satisfy the primal-dual optimality condition \eqref{eq:KKT}. Moreover, we have 
\begin{equation}\label{eq:ift}
 \frac{\partial y_\fd^*}{\partial x_\fd}= -J^{-1} \frac{\partial F}{\partial x_\fd}.
\end{equation}
Since
$$C_{\fd}^*(x_\fd, \tilde{\mu})= C_{\fd}(y_{\fd}^*(x_\fd); x_\fd, \tilde{\mu}),$$ 
we obtain
\begin{equation} \label{eq:partialDerivative}
\frac{\partial C^*_{\fd}}{\partial x_\fd}  = \frac{\partial C_{\fd}}{\partial y_\fd^*}^T \frac{\partial y_\fd^*}{\partial x_\fd},
\end{equation} 
where $\frac{\partial C_{\fd}}{\partial y_\fd^*}$ is the derivative of the subproblem cost function with respect to the local variable $y_\fd$.

We can also derive the Hessian of $C^*_\fd(x_{\fd};\tilde{\mu})$:
\begin{equation} \label{eq:partialHessian_vv}
\begin{split}
\frac{\partial^2 C^*_\fd}{\partial x_{\fd}^2}  & = \frac{\partial^2 y_\fd^*}{\partial x_{\fd}^2}^T\frac{\partial C_\fd}{\partial y_\fd^*} + \frac{\partial y_\fd^*}{\partial x_{\fd}}^T\frac{\partial^2 C_\fd}{\partial y_\fd^{*^2}}\frac{\partial y_\fd^*}{\partial x_{\fd}} \\
 & = (- \frac{\partial J^{-1}}{\partial x_{\fd}}\frac{\partial F}{\partial x_{\fd}} - J^{-1}\frac{\partial^2 F}{\partial x_{\fd}^2})\frac{\partial C_\fd}{\partial y_\fd^*} + \frac{\partial y_\fd^*}{\partial x_{\fd}}^T\frac{\partial^2 C_\fd}{\partial y_\fd^{*^2}}\frac{\partial y_\fd^*}{\partial x_{\fd}} \\
 & = -(-J^{-1}\frac{\partial J}{\partial x_{\fd}}J^{-1})\frac{\partial F}{\partial x_{\fd}}\frac{\partial C_\fd}{\partial y_\fd^*} + \frac{\partial y_\fd^*}{\partial x_{\fd}}^T\frac{\partial^2 C_\fd}{\partial y_\fd^{*^2}}\frac{\partial y_\fd^*}{\partial x_{\fd}}\\
 & = -J^{-1}\frac{\partial J}{\partial x_{\fd}}\frac{\partial y_\fd^*}{\partial x_{\fd}}\frac{\partial C_\fd}{\partial y_\fd^*} + \frac{\partial y_\fd^*}{\partial x_{\fd}}^T\frac{\partial^2 C_\fd}{\partial y_\fd^{*^2}}\frac{\partial y_\fd^*}{\partial x_{\fd}} .
\end{split}
\end{equation} 
The third step follows from the second step due to the identity $\frac{\partial J^{-1}}{\partial x_{\fd}} = -J^{-1}\frac{\partial J}{\partial x_{\fd}}J^{-1}$, and the term $J^{-1}\frac{\partial^2 F}{\partial x_\fd^2}$ is eliminated because $\frac{\partial^2 F}{\partial x_\fd^2} = 0$, since the variables in $x_\fd$, $V_{n_\fd}$ and $s_{n_\fd}$ appear only linearly in \eqref{eq:subproblem_original}.


An important observation is that the matrix $J$ in \eqref{eq:partialDerivative} and \eqref{eq:partialHessian_vv} is the same as the one used to compute Newton steps in the PDIPM. As a consequence, one can re-use the efficient implementation for solving \eqref{eq:Newton} that is already available in the PDIPM code, similar to the approach described in \cite{zavala2008fast}.
\vspace{-2pt}
\subsection{Two-Stage Algorithm}\label{eq:2stagealg}

The overall solution procedure is described in Algorithm~\ref{Alg:outer_alg}.  
\begin{algorithm}
  \caption{AC-OPF Two-stage Algorithm }\label{Alg:outer_alg}
  \textbf{Input:}  Initial smoothing parameter $\tilde\mu$, initial iterate $x_0$, stopping tolerance $\epsilon$.
  \begin{algorithmic}[1]
    \WHILE{$\tilde{\mu} > \epsilon$}
    \STATE\label{step:master} Solve master problem \eqref{eq:relaxed_master} with an NLP solver.

    Whenever NLP solver requires $C^*_\fd(x_\fd, \tilde{\mu})$, $\frac{\partial C^*_{\fd}}{\partial x_\fd}$, or $\frac{\partial^2 C^*_\fd}{\partial x_{\fd}^2}$, solve \eqref{eq:relaxed_subproblems} with Algorithm~\ref{Alg:inner_alg} and apply \eqref{eq:ift}, \eqref{eq:partialDerivative} and \eqref{eq:partialHessian_vv}.
      \STATE\label{step:outer_mu} Decrease smoothing parameter $\tilde\mu$.
      \ENDWHILE
  \end{algorithmic}
\end{algorithm}
It consists of solving a sequence of master problems \eqref{eq:relaxed_master} where the smoothing parameter $\tilde\mu$ is driven to zero.  An NLP method solves each of the master problems.  Whenever the NLP solver requires the value or derivatives of $C^*_\fd(x_\fd, \tilde{\mu})$, the values of $x_\fd$ corresponding to the current iterate are sent to the subproblems.  Each subproblem is then solved by a modified PDIPM (see Section~\ref{sec:PDIP}), and the derivatives are computed as described in Section~\ref{sec:2stagecom}.  These quantities are then sent back to the master problem NLP solver which continues its execution.

The optimal solutions of two consecutive master problems in Algorithm~\ref{Alg:outer_alg} can be expected to be close to each other, particularly when $\tilde\mu$ is small. Therefore, to aid convergence, the solution of a master problem is provided to the NLP solver as starting point for the solution of the next master problem after $\tilde\mu$ has been decreased. 
To best exploit this information, we prefer a sequential quadratic programming (SQP) method \cite{boggs1995sequential} over a nonlinear interior point method, since the latter is known to have inferior warm-starting capabilities \cite{potra2000interior}.
Particularly in the final stages of Algorithm~\ref{Alg:outer_alg} when $\tilde\mu$ changes only by a very small amount, we ideally want to encounter only very few iterations of the master problem solver. 
One of the conclusions of our numerical experiments is that this is indeed possible.
SQP methods generally require fewer function evaluations, an advantage in our situation where function evaluations are computationally expensive, and when an interior point QP solver is used to compute steps inside the SQP algorithm, large-scale instances can be solved.

\subsection{Infeasible Subproblems} \label{sec:penalty}
One difficulty of our two-stage decomposition approach is that subproblem \eqref{eq:subproblem} may become infeasible for a given $x_\fd$ during some master problem iteration.
As a remedy, 
using the notation from \eqref{eq:subproblem_original}, we let $\hat y_\fd= (V_{0_\fd}, s_{0_\fd})$ be a subvector of $y_\fd$  and introduce slack variables  $r$ and $t$ for the constraints \eqref{eq:constraint_v_equal}--\eqref{eq:constraint_s_equal} that force $\hat y_\fd$ to take the values $x_\fd$ prescribed by the master problem:
\begin{subequations} \label{eq:subproblem_l1}
\begin{align}
C^*_\fd(x_\fd) = \min_{y_\fd,r,t} & \quad  C_\fd(y_\fd;x_\fd) + \eta\;e^T(r+t)\\
s.t. & \quad \hat y_\fd - x_\fd = r - t, \label{eq:subproblem_l1_coupling}\\ 
& \quad g_\fd(y_\fd; x_\fd) = 0, \\
& \quad  h_\fd(y_\fd;x_\fd) \leq 0,  \\
& \quad r,t\geq 0.
\end{align}
\end{subequations}
Here, $e$ is the vector of all ones with appropriate dimension, and $\eta>0$ is a fixed parameter.
If the subproblem is feasible when $\hat y_\fd$ is not restricted, this problem is always feasible.
Clearly, at a (locally) optimal solution of \eqref{eq:subproblem_l1}, the new term in the objective measures the $\ell_1$-norm of the violation of the coupling constraints  \eqref{eq:constraint_v_equal}--\eqref{eq:constraint_s_equal}.
This is the standard $\ell_1$-norm penalty function formulation of a nonlinear optimization problem \cite{nocedal2006numerical}.
Penalty functions have been used in the past in two-stage decomposition approaches \cite{demiguel2006local,
braun1997development,alexandrov2002analytical}.
We choose the $\ell_1$-penalty function because it is ``exact'' in the sense that the optimal solution of  \eqref{eq:subproblem_l1} satisfies the original constraints \eqref{eq:constraint_v_equal}--\eqref{eq:constraint_s_equal} if the undecomposed problem is feasible and the penalty parameter $\eta$ is sufficiently large but finite \cite{nocedal2006numerical}.
The smoothing technique described in Section~\ref{sec:subproblem} is then applied to \eqref{eq:subproblem_l1} instead of \eqref{eq:subproblem}.

For the purpose of this paper we assume that a sufficiently large value for $\eta$ has been chosen.  A more comprehensive approach would include mechanisms that update the penalty parameter if it is too small \cite{byrd2010infeasibility}.

The penalty approach also makes it possible to extend our algorithm to settings with more than one coupling bus connecting the individual segments of the network obtained by decomposition.
In that case, coupling constraints are introduced for each coupling bus, together with slack variables as in \eqref{eq:subproblem_l1_coupling} and the corresponding penalty terms in the objective.  Again, when the original problem is feasible and the penalty parameter is large enough, the coupling constraints will be satisfied at an optimal solution computed by the decomposition algorithm, ensuring that the voltages and power flows of the connected network segments match.
In principle, this makes it possible to handle the decomposition of transmission systems into subnetworks.  We are planning to explore this possibility in future work.

\subsection{Convergence Guarantees}
The overall Algorithm~\ref{Alg:outer_alg} consists of three nested loops:

1) For given values of $x_\fd$ and $\tilde{\mu}$, convergence results for an appropriately chosen PDIPM method guarantee convergence of the algorithm solving the barrier problem \eqref{eq:relaxed_subproblems} to a local optimum, under standard assumptions that are typically satisfied by AC-OPF problems.  Recall that \eqref{eq:relaxed_subproblems} is feasible and hence the subproblem solution well-defined, due to the $\ell_1$-norm penalty formulation described in Section~\ref{sec:penalty} if we assume that the undecomposed problems is feasible.  

If the PDIPM always converges to a unique global minimum of \eqref{eq:relaxed_subproblems}, then $C^*_\fd(x_\fd,\tilde{\mu})$ is uniquely defined and, as discussed in Section~\ref{sec:2stagecom}, differentiable.
This assumption is reasonable since it has been observed that, in practice, a PDIPM applied to an AC-OPF problem usually converges to the global minimizer
\cite{7763860}.

2) For a fixed value of $\tilde{\mu}$, the modified master problem \eqref{eq:relaxed_master} is a nonlinear optimization problem with differentiable problem functions. An appropriately chosen NLP solver, such as an SQP method, will converge to a local optimum (if the original two-state problem is feasible).

3) Finally, to understand the convergence of the overall Algorithm~\ref{Alg:outer_alg}, we cite a result from \cite{demiguel2006local} that discusses the existence of local solutions of the two-stage problem:

\begin{thm} Let $(x^*, \{x_\fd^*\}_{\fd \in \fD}, \{y_\fd^*\}_{\fd \in \fD})$ be a minimizer of the undecomposed problem satisfying the nondegeneracy conditions in \cite{demiguel2006local}. Then there exists a locally unique trajectory $(x^*(\tilde{\mu}), x^*_\fd(\tilde{\mu}))$ of minimizers to \eqref{eq:relaxed_master}, such
$~\lim_{\tilde{\mu} \to 0} x^*(\tilde{\mu}) = x^*$ and $\lim_{\tilde{\mu} \to 0} x^*_\fd(\tilde{\mu}) = x_\fd^*.$
\end{thm}

Therefore, Algorithm~\ref{Alg:outer_alg} converges to a  local minimizer $(x^*, \{x_\fd^*\}_{\fd \in \fD})$ if the master problem solver eventually returns optimal solutions corresponding to its local trajectory.
Again, since local solvers typically find the global solutions of AC-OPF problems, this is a reasonable assumption.

\section{Numerical Example} \label{sec:numerical}
\subsection{Implementation Details}
Our algorithm is implemented in MATLAB R2020a. The SQP algorithm in the optimization package KNITRO 10.3.0 \cite{byrd2006knitro} is employed to solve the master problem. 
When $\tilde{\mu}$ decreases, we initialize the master problem with the optimal solution from the previous $\tilde{\mu}$. 
To solve the subproblems, we implemented a basic version of the PDIPM Algorithm~\ref{Alg:inner_alg} given in Section \ref{sec:subproblem}. Warm-start initialization is also used for the solution of subproblems \eqref{eq:relaxed_subproblems}, where at each NLP iteration, we initialize the PDIPM with the primal-dual optimal solutions from the last NLP iteration. Experiments show that with this warm-start initialization, the PDIPM algorithm always converges within 10 iterations with a tolerance  $10^{-12}$. The penalty parameter was chosen as $\eta=9$, and the tolerance of the SQP method was kept at its default value $10^{-6}$.

The experiments were conducted on a Linux workstation with two 20-core Intel Xeon processors, running at 2.40 GHz, and with 256GB RAM.  We had exclusive access to this machine so that the time measurements are not affected by other processes. For all experiments we report the wall clock time averaged over three runs with identical data in order to reduce the effect of random variations in run times.  Setup time, such as reading and generating data, is not included. 

\subsection{Numerical Results}

This section describes numerical experiments obtained with a combined transmission and distribution network model. Since the focus of this paper is the decomposition algorithm, the distribution model is simplified as a single-phase model. Our approach can be 
 extended to three-phase distribution system in a straight-forward manner, by changing the communication variable $x_\fd$ and the subproblem formulation. 


The test system consists of one transmission system and $|\fD|$ distribution systems. We took the IEEE-118 system \cite{christie1993ieee} as the transmission system, and attached a distribution system at each of the 64 buses that have degree of at most 2.  As distribution system we chose a meshed 143-bus system based on the 135-bus system in \cite{hijazi2015optimal, guimaraes2005reconfiguration} with 8 generators and 156 branches.
The resulting combined network then has $118 + 64\times (143-1)=9,206$ buses.
%
The data in the distribution systems is scaled to match the loads of the corresponding buses in the transmission system. The solutions of our decomposition algorithm on this test data recover the solutions computed by MATPOWER \cite{zimmerman1997matpower} for the undecomposed system.

To stress our algorithm with very large instances and permit scaling experiments, we consider a stochastic formulation that accounts for uncertainties in the net demand in the distribution systems.
We generated $N$ scenarios for each subnetwork $\fd \in \fD$ with $\pm 10 \%$ fluctuations in subnetwork loads. 
%
More specifically, for scenario $i$, we chose a vector of fluctuations $\xi_{\fd,i}$ randomly from a uniform distribution and multiply all loads in the subnetwork $\fd$ by the corresponding components of $\xi_{\fd,i}$.
The goal is the minimization of the expected total generation cost, estimated by sample average approximation. In this case, the objective function of the two-stage problem \eqref{eq:relaxed_master} becomes
\begin{equation}
C(x) + \sum_{\fd \in \fD}\frac{1}{N}\sum_{i=1}^N C^*_\fd(x_\fd, \tilde{\mu}; \xi_{\fd,i})
\end{equation}
with $|\fD| \times N$ many subproblems, each one corresponding to a particular scenario for a given subnetwork.
Here, $C^*_\fd(x_\fd, \tilde{\mu}; \xi_{\fd,i})$ denotes the subproblem objective \eqref{eq:relaxed_subproblems_obj}, now for the vector of demand fluctuations $\xi_{\fd,i}$.

In this model, the first-stage decision variables correspond to the transmission network and are fixed for all potential second-stage scenarios.  Only the second-stage variables, which correspond to the distribution networks, can be adjusted once the uncertainty is realized.  In other words, we assume that the distribution systems react in a way so that the transmission system is not affected.  This corresponds to a setting where an active distribution network strives to meet the pre-negotiated demand profile in order to avoid penalty payments.

\subsubsection{Performance with increasing number of scenarios}
\begin{figure}[!t]
  \centering
  \includegraphics[width=0.45\textwidth, height=0.32\textwidth]{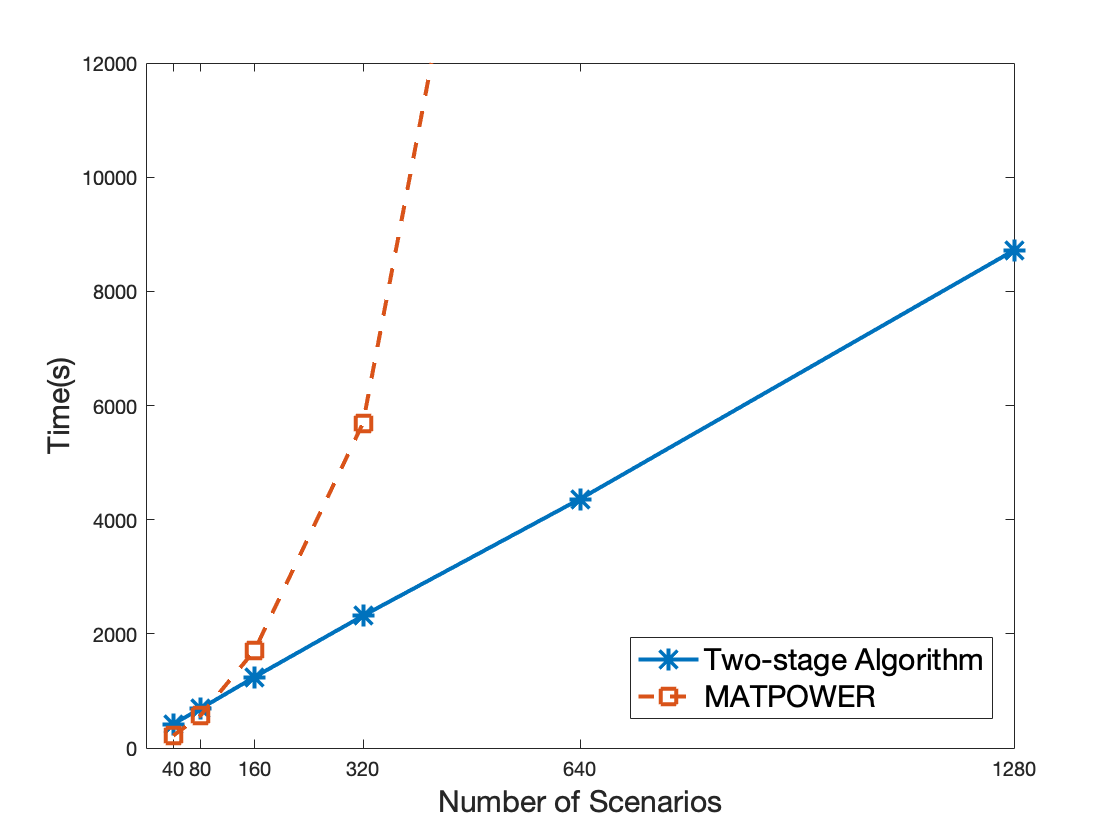}
  \caption{Computation time as the number of scenarios is increased.}
  \label{fig:timeVSscenario}
\end{figure}

To explore the scalability of the two-stage decomposition approach, we ran our implementation of Algorithm~\ref{Alg:outer_alg} on the integrated power system with  $|\fD| = 64$  subnetworks. The number of scenarios per subnetwork ranged from 40 to 1,280. These experiments were performed in a parallel computing setting with 40 threads, and the wall clock times required to converge are measured.  The largest data set results in an equivalent undecomposed network with
$118 + 64 \times 1280 \times (143-1) = 11,632,758$ buses in total. We repeated the same experiment using MATPOWER and compare the computation time in Figure~\ref{fig:timeVSscenario}. Our two-stage algorithm and MATPOWER converged to the same objective values on these test cases.
We observe that the computation time of the two-stage decomposition algorithm increases linearly in number of scenarios, as expected, while it increases superlinearly for MATPOWER. In particular, MATPOWER took 7.9 hours for $N=640$ and failed to converge for the largest instance within a time limit of 24 hours.  In contrast, the decomposition algorithm is able to solve the largest instance in about 2.3 hours. The experiment indicates that our approach scales better than a centralized solution as the number of subproblems increases, and it succeeds in solving large-scale problems that can be out of reach for an undecomposed formulation.

In these experiments, the decomposition algorithm was given a flat start as starting point, which is less favorable than the starting point included in the MATPOWER models that was used by the MATPOWER optimization algorithm.  Initializing MATPOWER with a flat start would have led to a larger number of iterations and longer running times.

\begin{figure}[!t]
  \centering
  \includegraphics[width=0.45\textwidth, height=0.32\textwidth]{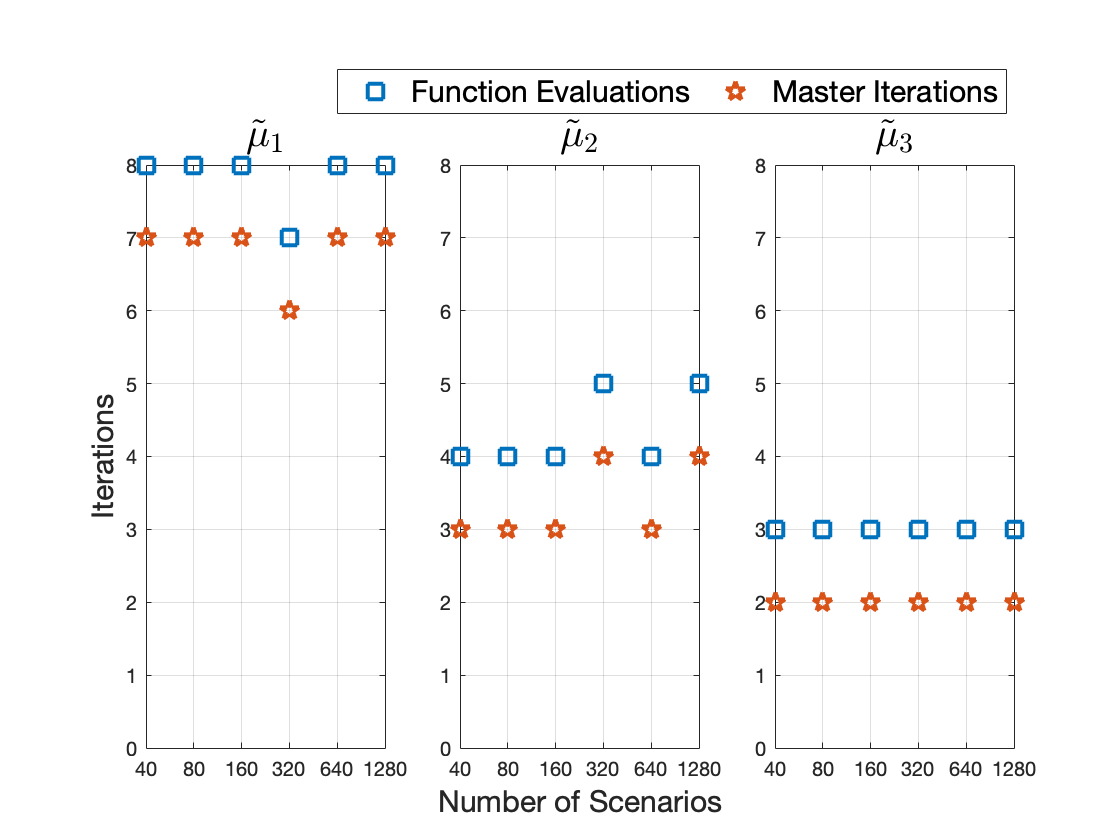}
  \caption{Number of master problem iterations and function evaluations for different problem sizes and smoothing parameter values.}
  \label{fig:scatter_scenario}
\end{figure}

Moreover, as the number of scenarios grows, the number of master problem iterations remains constant. In our implementation, we set the initial value of $\tilde{\mu}_1$ to be $10^{-2}$, and then sequentially decrease its value in Step~\ref{s:muupdate} of Algorithm~\ref{Alg:outer_alg} to $\tilde\mu_2=10^{-3}$ and $\tilde\mu_3=10^{-6}$, where the last value corresponds to the final tolerance.
For each $\tilde{\mu}_i$, we observed that the master problem SQP solver \eqref{eq:relaxed_master} converges locally at a superlinear rate.  This can be expected because second-order derivatives are provided \cite{nocedal2006numerical}.
Figure \ref{fig:scatter_scenario} shows the number of SQP iterations and the number of function evaluations for the different numbers of scenarios. 
We find that the number of iterations is small, at most 7, and that only one function evaluation was required per iteration after the evaluation at the starting point.  This means that the line search procedure in the SQP algorithm accepted every trial point.
We also observe that the warm start initialization described in Section~\ref{eq:2stagealg} is highly efficient.
In particular, for the final smoothing parameter $\tilde\mu_3$, the master problem converged within just 2 iterations.

\subsubsection{Performance with increasing number of threads}
\begin{figure}[!t]
  \centering
  \includegraphics[width=0.45\textwidth, height=0.32\textwidth]{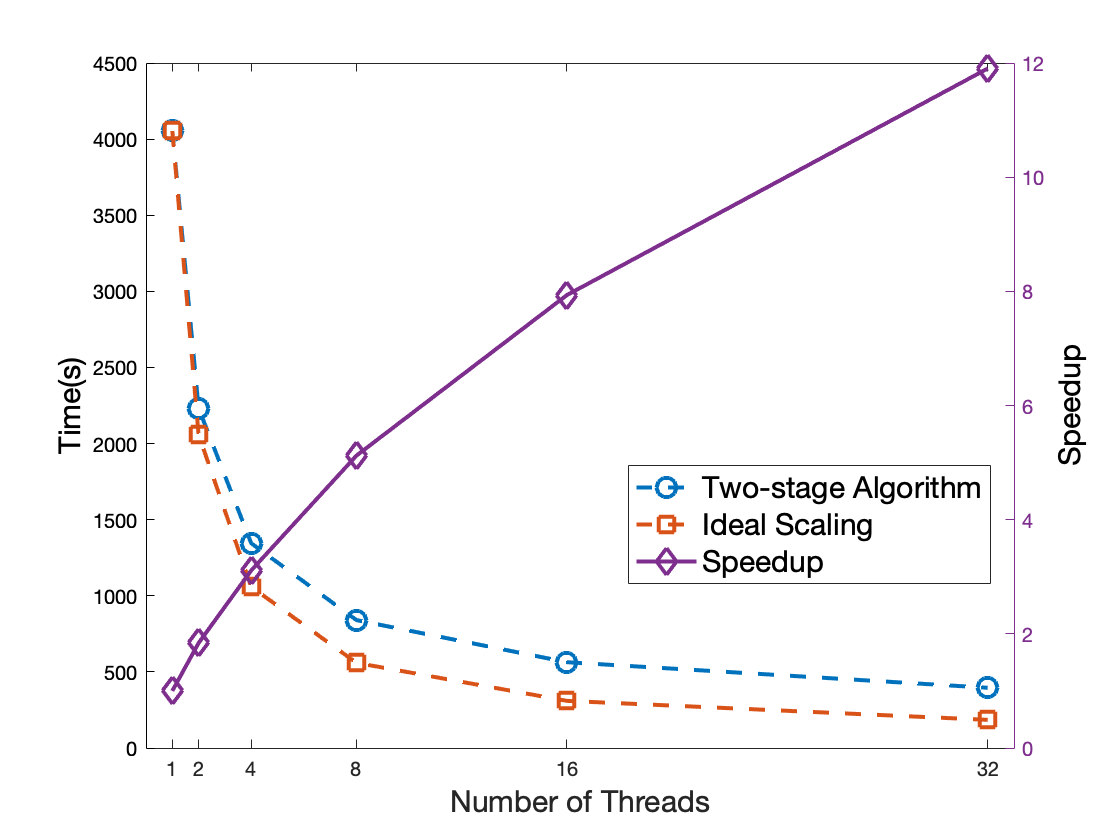}
  \caption{Strong scaling test: Fixed number of subnetworks and scenarios.}
  \label{fig:timeVScore}
\end{figure}

To explore the speedup obtained by parallel computations, we performed both strong and weak scaling tests for the integrated system. The strong scaling test is conducted using a fixed test case with $|\fD| = 64$ subnetworks and $N=64$ scenarios for each subnetwork, i.e., with a total of about 553k buses for this experiment, as we increase the number of parallel threads. Figure \ref{fig:timeVScore} plots the wall clock time against the number of threads we used. 
We see that the computation time is significantly reduced by solving subproblems in parallel.
However, in our experiments, the speedup was less than ideal.  Using 32 threads lead to a reduction by a factor of only 12.  We believe that this can only partially be attributed to the unparallelized portion of the code for the solution of the master problem (KNITRO took a total of around 40 secs), and to the usual degradation due to the competition of the cores for the shared hardware resources, such as cache and the memory bus. 
The major cause appears to be overhead in MATLAB's parallelization infrastructure.

\begin{figure}[!t]
  \centering
  \includegraphics[width=0.45\textwidth, height=0.32\textwidth]{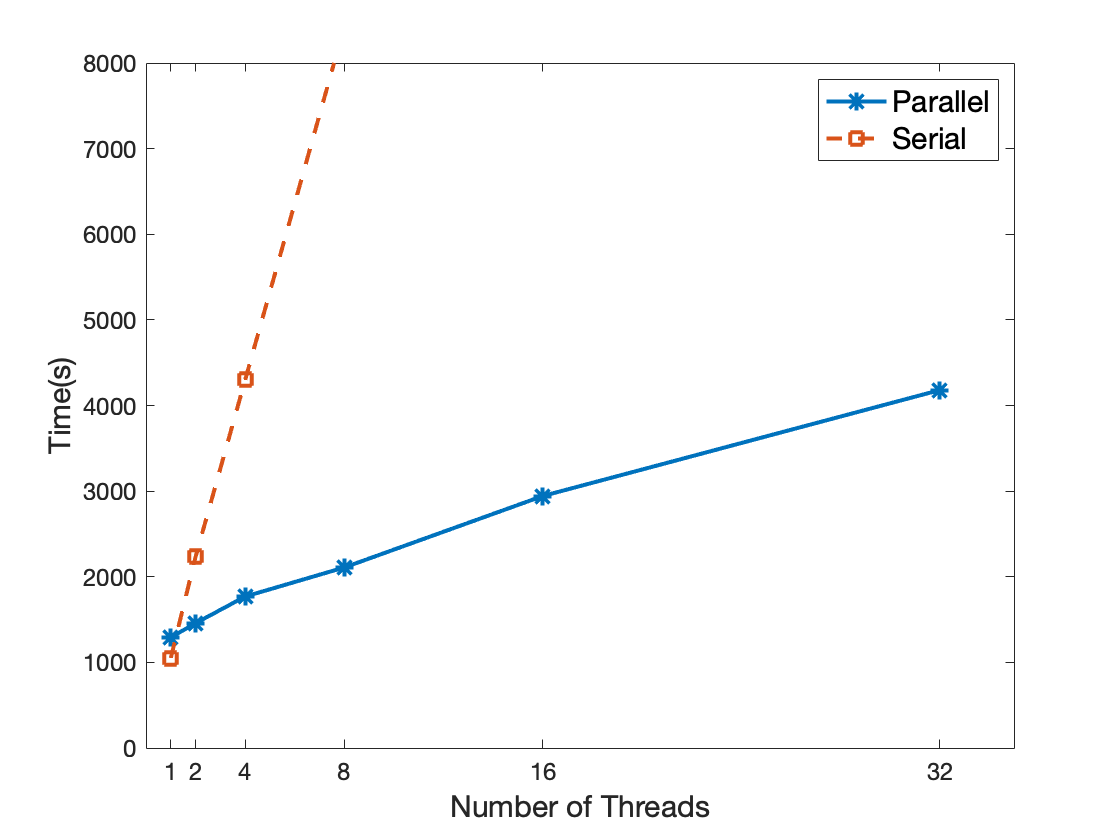}
  \caption{Weak scaling test: Increasing number of scenarios.}
  \label{fig:ratio1}
\end{figure}

\begin{figure}[!t]
  \centering
  \includegraphics[width=0.45\textwidth, height=0.32\textwidth]{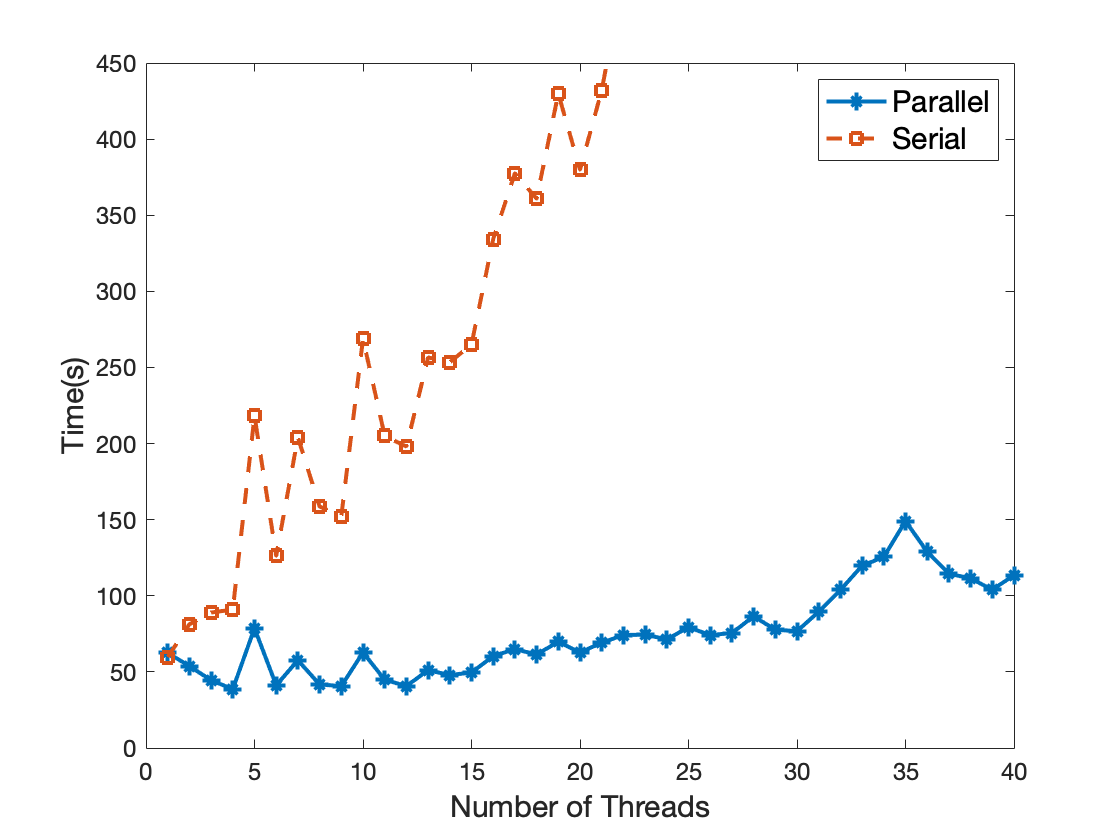}
  \caption{Weak scaling test: Increasing number of subnetworks.}
  \label{fig:ratio2}
\end{figure}

For the weak scaling test, two different settings are considered. The first scaling test fixes the number of $|\fD| = 64$ subnetworks and increases the number of $N = 20\times \text{(number of threads)}$ scenarios.  The second test considers the case with a fixed number of $N = 20$ scenarios and an increasing number of $|\fD| = 1 \times \text{(number of threads)}$ subnetworks. The wall clock times are shown in
Figures~\ref{fig:ratio1} and \ref{fig:ratio2}.  For comparison, we also present the run times when a single thread is used.
With ideal scaling, the time for the parallel runs should remain constant, independent of the number of threads.
For the experiment depicted in Figure~\ref{fig:ratio1}, we see a deviation from the ideal speedup by a factor of about 3, similar to what is observed in Figure~\ref{fig:timeVScore}.
On the other hand, close to ideal scaling is observed in Figure~\ref{fig:ratio2}. Note that the computation time in Figure \ref{fig:ratio2} is not monotone in the number of threads.  This is due to the fact that the network topology changes as number of subnetworks attached to the master network increases and the number of master problem iterations varies.
In both cases we see that, as expected, the time for the serial run increases at least linearly with the size of the undecomposed equivalents.

\section{Summary} \label{sec:conclusion}
This paper proposes a novel 
two-stage optimization algorithm
that was applied to partitions of a power network into a master transmission network and a set of distribution subnetworks.  Stochastic instances were handled by replicating subnetworks with different realizations of the uncertain parameters.

By introducing a barrier term  as a smoothing technique for the subproblems, the response of the second-stage value function becomes differentiable with respect to the master problem variables.  As a consequence, efficient and reliable existing nonlinear optimization packages with fast local convergence properties can be utilized for the master problem.
The exploitation of warm-start capabilities of an SQP solver significantly accelerates the solution of subsequent master problems in which the value of the smoothing parameter is decreased. 

The second-stage AC-OPF problems can be solved 
with primal-dual interior point methods which exhibit fast local convergence guarantees as well.  Existing algorithm implementations can be used after minor modifications of their termination criteria.  First- and second-order derivatives of the subproblem response with respect to the master problem variables can be derived via the implicit function theorem and computed efficiently using the Jacobian matrix of the primal-dual optimality conditions, which is already constructed and factorized within the interior point solver.

The framework is naturally able to exploit parallel computing resources by distributing subproblem calculations. Only a small amount of communication is required between the master problem and subproblems.
Our experiments show that the approach is able to solve large-scale instances with more than 11 million buses, which is out of reach for centralized solutions, such as MATPOWER.  We also showed that it scales well with an increase of the number of parallel computing threads.

The proposed algorithm is not limited to the specific optimization of TSO-DSO interactions.  It is designed for general two-stage optimization problems with continuous variables and is likely to offer an efficient solution framework for many other applications. 

One direction of future work is towards a more realistic setting in which TSO and DSO do not have the same objectives with an extension of the smoothing technique to bilevel formulations. 

\section*{Acknowledgements}

We thank Hassan Hijazi for providing us with the model of the 143-bus distribution system.  We also thank Dan Molzahn for his comments on an early draft of this article.





\ifCLASSOPTIONcaptionsoff
  \newpage
\fi




\bibliographystyle{IEEEtran_bib}
\bibliography{IEEEabrv,ATwoStageDecompositionApproachforACOptimalPowerFlow_Tu}
%

%

\vspace*{-4.5\baselineskip}
\begin{IEEEbiography}[{\includegraphics[width=1in,height=1.25in,clip,keepaspectratio]{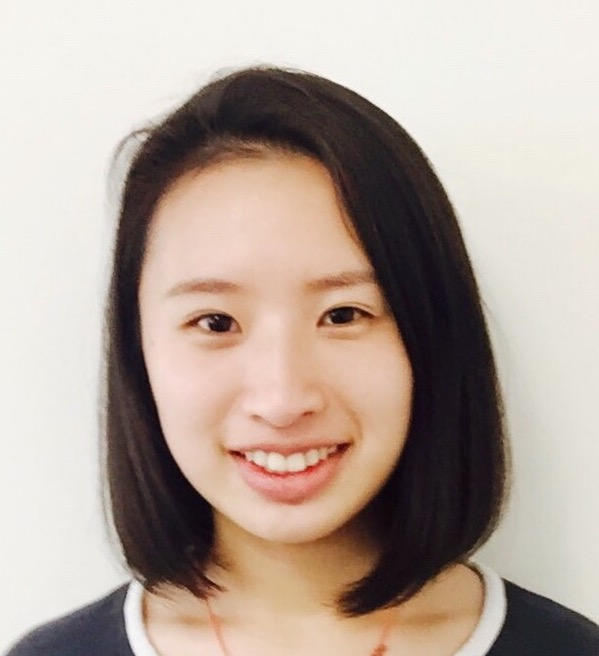}}]{Shenyinying Tu}
received the B.S. degree in Applied Mathematics from University of California, Los Angeles, U.S., in 2016. She is
currently pursuing the Ph.D. degree in Industrial Engineering ad Management Science, Northwestern University. Her research
interests include optimal power flow, nonlinear optimization.z
\end{IEEEbiography}

\vspace*{-4\baselineskip}
\begin{IEEEbiography}[{\includegraphics[width=1in,height=1.25in,clip,keepaspectratio]{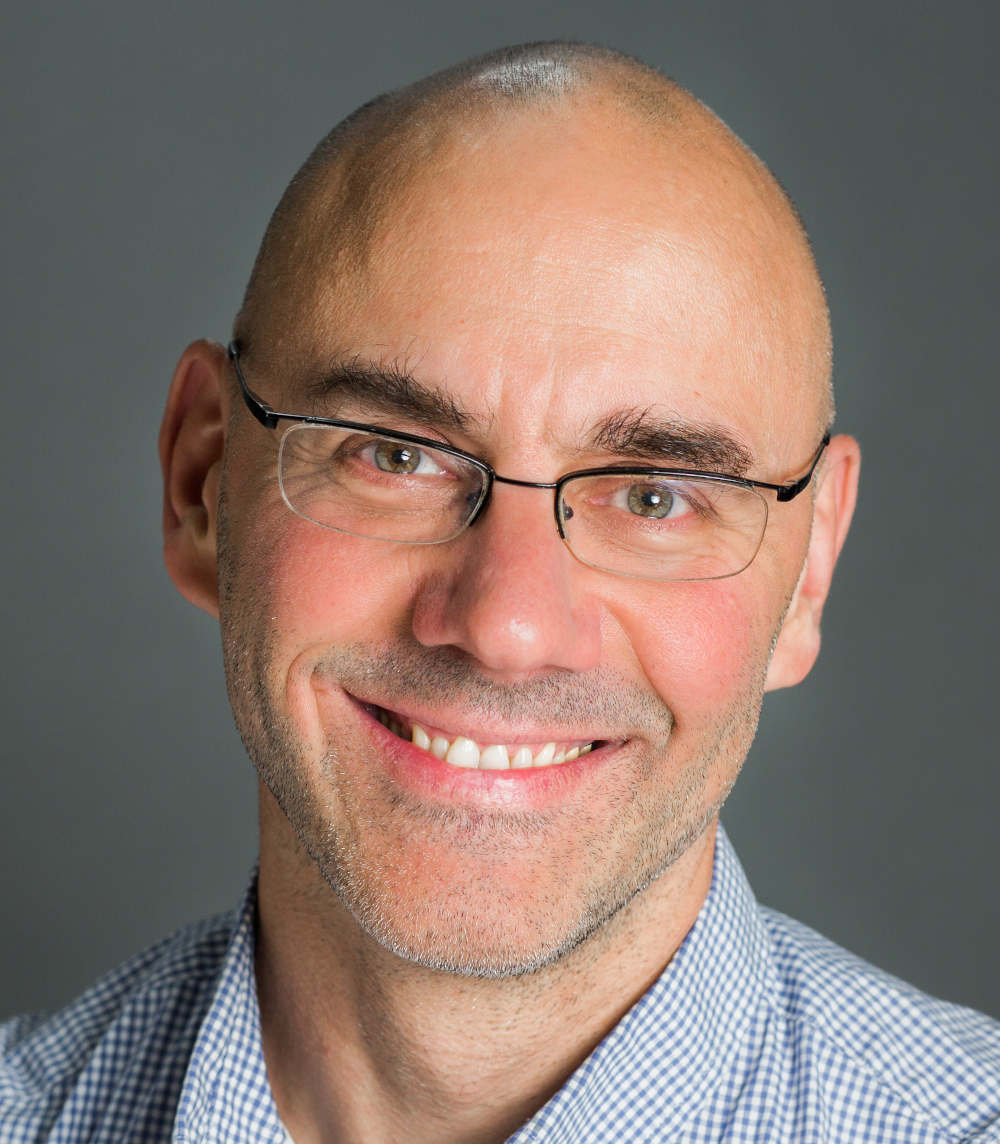}}]{Andreas W\"achter}
is a Professor at the Industrial Engineering and Management Sciences Department at Northwestern University.  He obtained his master's degree in Mathematics at the University of Cologne, Germany, in 1997, and this Ph.D. in Chemical Engineering at Carnegie Mellon University in 2002.  Before joining Northwestern University in 2011, he was a Research Staff Member in the Department of Mathematical Sciences at IBM Research in Yorktown Heights, NY.  His research interests include the design, analysis, implementation and application of numerical algorithms for nonlinear continuous and mixed-integer optimization.  He is a recipient of the 2011 Wilkinson Prize for Numerical Software and the 2009 Informs Computing Society Prize for his work on the open-source optimization package Ipopt.  He is currently spending a year at Los Alamos National Laboratory as the Ulam Fellow.
\end{IEEEbiography}

\vspace*{-4\baselineskip}

\begin{IEEEbiography}[{\includegraphics[width=1in,height=1.25in,clip,keepaspectratio]{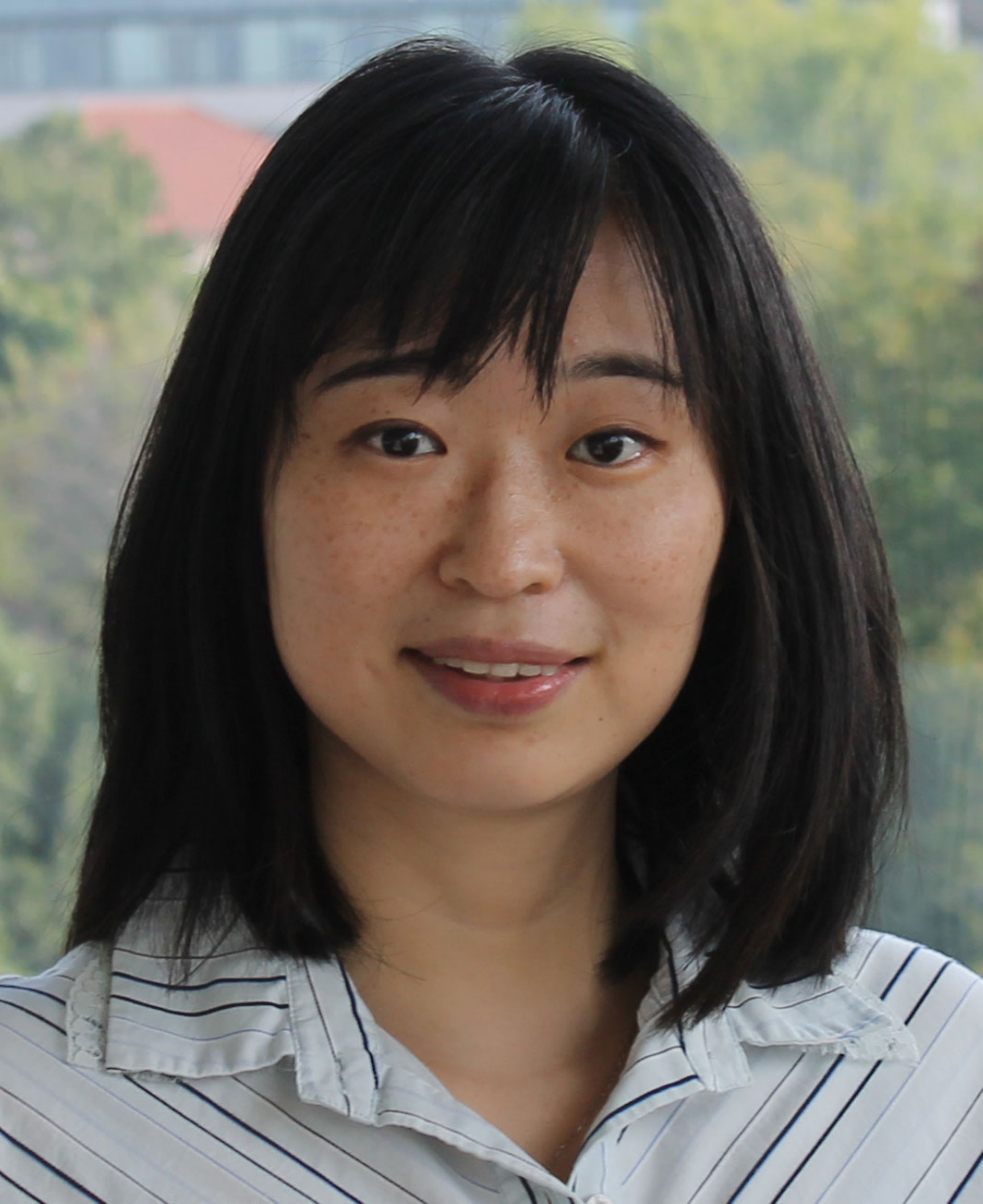}}]{Ermin Wei}
is currently an Assistant Professor at the Electrical and Computer Engineering Department and Industrial Engineering and Management Sciences Department of Northwestern University. She completed her PhD studies in Electrical Engineering and Computer Science at MIT in 2014, advised by Professor Asu Ozdaglar, where she also obtained her M.S.. She received her undergraduate triple degree in Computer Engineering, Finance and Mathematics with a minor in German, from University of Maryland, College Park. Wei has received many awards, including the Graduate Women of Excellence Award, second place prize in Ernst A. Guillemen Thesis Award and Alpha Lambda Delta National Academic Honor Society Betty Jo Budson Fellowship. Wei's research interests include distributed optimization methods, convex optimization and analysis, smart grid, communication systems and energy networks and market economic analysis.
\end{IEEEbiography}




\end{document}